\documentclass{amsart}
\usepackage[latin1]{inputenc}
\usepackage{amsmath}
\usepackage{amssymb}

\usepackage{hyperref}

\begin{document}

\title{Making research on symmetric functions with MuPAD-Combinat}

\author{Francois Descouens}

\maketitle

\begin{abstract}
  We report on the 2005 AIM workshop ``Generalized Kostka
  Polynomials``, which gathered 20 researchers in the active area of
  $q,t$-analogues of symmetric functions. Our goal is to present a
  typical use-case of the open source package MuPAD-Combinat in a
  research environment.
\end{abstract}

\section{Introduction}
The project \texttt{MuPAD-Combinat} (see
\url{http://mupad-combinat.sf.net/} and \cite{HTm}), born in spring
2001 under the leadership of F. Hivert and N.  Thi\'ery, is an open
source package for making algebraic combinatorics using the computer
algebra system MuPAD (see \url{www.mupad.de} and \cite{Fuch} for more
details). The main goal of this package is to
bring an open source flexible and easily extensible toolbox for
checking conjectures in a short programming time.  This package
contains a large collection of tools as implementations of classical
combinatorial objects (partitions, Young tableaux, trees, $\ldots$),
computations in combinatorial Hopf algebras (in particular symmetric
functions), manipulations of graphs, automata, $\ldots$

The study of symmetric functions is a major historical field in
algebraic combinatorics (\cite{Mcd}) and some people are mainly
interested in generalizations of Kostka polynomials. In section 2, we
give briefly mathematical definitions about different objects we work
on. We illustrate how a computer algebra system which contains
classical combinatorial objects, evolutive implementation of symmetric
functions and easy way to incorporate program written in C++ can be
used for our research.

In section 3, we explain the design of our implementation of symmetric
functions and some of our technical choices. We also give some
examples of computations and show some advanced functionalities as
adding new bases on the fly (using different characterizations) or
implementing new operators on symmetric functions. Finally, in section
4 we describe how we incorporate programs into a coherent design using
MAPITL library. The technical concepts used in section 3 and 4 are
essentially classical, but we want to stress on their integration into
the package in order to have a powerfull tool for making efficient.
The difficulty is to find the appropriate combination that yields an
intuitive yet flexible and powerful research tool. We show that our
choices are promising in a real situation of collaboration at a
workshop on generalized Kostka polynomials in Palo Alto, California
organized by the American Institute of Mathematics in July 2005 (see
\url{ http://www.aimath.org/WWN/kostka/ } for more informations on
this event).

\section{Symmetric functions and Kostka polynomials}
\subsubsection{Basic definitions}
A symmetric polynomial in variables $X=\{x_1,\ldots,x_n\}$ is a
polynomial in $X$ invariant under permutations of variables. When is
infinite, we call symmetric function such a polynomial. The set of
symmetric functions with coefficients in $\mathbb{C}(t)$, denoted
$\Lambda_{t}$, is a graded algebra with respect to the degree of
polynomials, i.e
$$
\Lambda_{t}=\oplus_{n\ge 0}\Lambda_{t}^n.$$
For all $n\ge 0$, the
dimension of $\Lambda_{t}^n$ is the number of partitions of $n$ (a
partition of a positive integer $n$, written $\lambda\vdash n$, is a
decreasing sequence of positive integers with sum $n$). One main basis
of this algebra is constituted with monomial functions defined for all
partitions by $\lambda=(\lambda_1,\ldots,\lambda_n)$, by
$$m_\lambda=\sum_{v\in O(\lambda)}x_1^{v_1}\ldots x_n^{v_n},$$
where
$O(\lambda)$ represents the set of all the permutations of $\lambda$.
Another interesting basis consists of the symmetric powersums,
defined for all partitions $\lambda$ by
$$p_\lambda(X)=p_{\lambda_1}\ldots p_{\lambda_n}\quad
\text{where}\quad p_{\lambda_i}=\sum_{j=1}^n x_j^{\lambda_i}.$$
A
scalar product on $\Lambda_{t}$ can be uniquely defined by
$$
(p_\lambda, p_\mu)=\delta_{\lambda,\mu}\prod_{i\ge 1}(m_i)!\ 
i^{m_i(\lambda)},$$
where $m_i(\lambda)$ represents the multiplicity
of part $i$ in partition $\lambda$. Applying Gram-Schmidt process of
orthonormalization on the monomial basis with respect to this scalar
product yields us the basis of Schur functions $(s_\lambda)_\lambda$.
These functions are intensively studied (from an algebraic and
combinatorial point of view) and one of the beautiful results is the
Littlewood-Richardson rule, a combinatorial interpretation of the
product of two Schur functions.
\subsubsection{Kostka polynomials}
One can introduce a $t$-deformation of the previous scalar product by
$$(p_\lambda, p_\mu)_t=\delta_{\lambda,\mu}\prod_{i\ge 1}(m_i)!\ 
i^{m_i(\lambda)}\prod_{i=1}^{l(\lambda)}\frac{1}{1-t^{\lambda_i}}.
$$
The orthogonalization of the monomial basis with respect to this
scalar product defines the Hall-Littlewood functions $P_\lambda(X;t)$.
There exist two other families of Hall-Littlewood functions: on the
one hand, the $Q_\lambda(X;t)$ which are the dual elements of
$P_\lambda(X;t)$ with respect to the scalar product $(\, ,\,)_t$ and
on the other hand $Q^{'}_\lambda(X;t)$ which are the dual elements of
$P_\lambda(X;t)$ with respect to $(\, ,\,)$. The expansion of the
$Q^{'}_\lambda$ on Schur functions is an algebraic way to define the
Kostka polynomials $K_{\lambda,\mu}(t)$
$$Q^{'}_\lambda(X;t)=\sum_{\mu\vdash \vert \lambda \vert }K_{\lambda,\mu}(t)s_\mu .$$
These
polynomials have also a combinatorial interpretation, using the charge
on semi-standard Young tableaux in \cite{LS}, or the rigged configurations
introduced by Kerov, Kirillov and Reshetikhin in \cite{KKR}.
\subsubsection{Generalizations of Kostka polynomials} 
Using $k$-ribbon tableaux introduced by Lascoux, Leclerc and Thibon in
\cite{LLT}, we define for each positive integer $k$ and each partition
a particular symmetric functions $H_\lambda^{(k)}(X;t)$.  Their
expansion on Schur functions gives us a way to define an increasing
filtration of Kostka polynomials $K_{\lambda,\mu}^{(k)}(t)$
$$
H_\lambda^{(k)}(X;t)=\sum_{\mu\vdash \vert \lambda \vert} K_{\lambda,\mu}^{(k)}(t)s_\mu .$$
The expansion of $K_{\lambda,\mu}^{(k)}(t)$ on the monomial basis is a
way to define unrestricted generalized Kostka polynomials. In order to
generalize Kostka polynomials, there exist other combinatorial ways
(with unrestricted rigged configurations recently introduced by L.
Deka and A. Schilling in \cite{Schilling3}) and algebraic ways (using
the theory of crystal bases for quantum groups of type $A_n$) to
generalize Kostka polynomials. Other generalizations are given by M.
Zabrocki using creation operators in \cite{Zabrocki} and by L.
Lapointe and J.  Morse introducing a $t$-deformation of $k$-Schur
functions in \cite{Morse}. An interesting problem, explained in
\cite{Schilling1}, is to show that these generalizations coincide in
some particular cases.

\section{Implementation of symmetric functions}
\subsection{Design goals}

Symmetric functions can be represented in many different ways, and in
particular in different basis (powersum,elementary, Schur,
monomials,$\ldots$). As usual in computer science, it is essential at
each point to use the appropriate representation, both for efficiency
and interpretation of the results. What make the situation specific is
the number of those representations. In particular, it is neither
practical nor sometimes possible to implement explicitely all
conversions. Instead we want to be able to only implement a few and
deduce the others by compositions or linear tranformations such as
inversion or transposition. Furthermore to thratend, the user and also
programmers should not need to know which conversion are implemented
because this information is too volatile.  For example, the addition
of symmetric function not given in the same basis is possible

\scriptsize 
\begin{verbatim}
>> S::s([2,1]) + S::QP([2,1]) + S::p([2,1]);

               (t + 4) m[1, 1, 1] + (t + 1) m[3] + (t + 3) m[2, 1]
\end{verbatim}
\normalsize The monomial basis has been choosen by the system because
it minimizes the number conversions (the cost of the conversion is not
taken into account). The same process accurs for making conversion
between two bases as the expansion of the Hall-Littlewood function
$Q^{'}_{311}$ on the monomial basis

\scriptsize
\begin{verbatim}
  >> S::m(S::QP([2,1]));

               (t + 2) m[1, 1, 1] + (t + 1) m[2, 1] + t m[3]
\end{verbatim}
\normalsize We also consider operators on symmetric functions. In most
cases, they are easier to define (and consequently to implement) on
one of the bases but not on all of them. Consequently, the system uses
implicit conversions between bases in order to apply operators on any
given basis. The main technologies used are linear algebra and
overloading mechanisms which are not essentially new, but a great
effort is made in order to make manipulation intuitive.

\subsection{General overview of the implementation}
In our design, for each basis of symmetric functions there is a domain
(in the category 
\texttt{Cat::GradedHopfAlgebraWithBasis}) which represents
the space of symmetric functions expanded on this basis. For example,
here is the example of implementation of the complete basis
\scriptsize
\begin{verbatim}
domain SymT::complete(R: DOM_DOMAIN)
    inherits SymT::common(R);
    category Cat::GradedHopfAlgebraWithBasis(R), Cat::CommutativeRing;
    info_str := "Domain for symmetric functions expanded on complete basis";

    basisName := hold(h);

    // Implementation of multiplication (complete basis is a multiplicative basis)
    mult2Basis := dom::term@revert@sort@_concat; 
    ...
end_domain:
\end{verbatim}
\normalsize

The domain \texttt{Sym} of symmetric functions, representing symmetric
functions in whatever representation, is in category
\texttt{Cat::HopfAlgebraWith\- SeveralBases}. This category helps us
in defining implicit conversions between bases

\scriptsize
\begin{verbatim}
domain Sym(R=Dom::ExpressionField()) 
    inherits Dom::BaseDomain;
    category Cat::HopfAlgebraWithSeveralBases(R), Cat::CommutativeRing;
    ...

    // Each domain corresponding to a basis is declared 
    h := SymT::complete(dom::coeffRing, Options);

    // Implementation of different conversions between bases (stored in a table)

    basisChangesBasis := 
    table(
          ...
          // Explicit combinatorial conversions 
          (dom::QP, dom::s)   = (part ->(_plus(combinat::tableaux::kostkaPol(mu, part, dom::vHL)
                                               * dom::s(mu) $ mu in
                                                 combinat::partitions::list(_plus(op(part)))))),
          ...
          (dom::McdP, dom::m) = (part ->((dom::GramSchmidt(dom::m, _plus(op(part)), 
                                          dom::scalartq))[op(part)])),
          ...
          // Dual conversions and inverse conversions
          (dom::s, dom::QP) = dom::invertBasisChange(dom::QP, dom::s),
          (dom::s, dom::m) = dom::transposeBasisChange(dom::h, dom::s, dom::s, dom::m),
         )
\end{verbatim}
\normalsize Note that only some conversions are implemented in this
table; the overloading mechanism is in charge of finding the shortest
number of intermediate conversions needed (it doesn't take into
account the cost of each conversion). In order to use
$(q,t)$-deformation of symmetric function, we can declare \scriptsize
\begin{verbatim}
>> S:=examples::SymmetricFunctions(Dom::ExpressionFieldWithDegreeOneElements([t,q]), 
                                                                             vHL=t, vMcd=q);
\end{verbatim}
\normalsize
\subsubsection{Adding new bases on the fly} 
In order to add new basis on the fly, we define a generic domain.
\scriptsize
\begin{verbatim}
  domain SymT::NewBasis(R: DOM_DOMAIN, DomName: DOM_STRING)
      inherits SymT::common(R);
      category Cat::GradedHopfAlgebraWithBasis(R), Cat::CommutativeRing;
      info_str := "Domain for symmetric functions expanded on a new added basis";
        
      basisName := text2expr(DomName);  
        
  end_domain:
\end{verbatim}
\normalsize
In order to add a new basis named $E$ for example, we can use
\scriptsize
\begin{verbatim}
  >> B := S::newBasis(S::coeffRing, ``E'');
\end{verbatim}
\normalsize and we define the change of bases we want. Let suppose
that the change of basis between B and monomial basis corresponds to
the function  \texttt{testChange}  
\scriptsize
\begin{verbatim}
  >> S::declareBasisChangeBasis((B, dom::m, testChange);

  >> S::declareBasisChangeBasis(dom::m, B, dom::invertBasisChange(B, dom::m));
\end{verbatim}
\normalsize

\subsubsection{Adding new operators}
If we want to define new operators on symmetric functions, we only
have to define their action on a particular basis, of course for
efficiency. First, we declare our operators as an overloaded operator
\scriptsize
\begin{verbatim}
   >> newOp := operators::overloaded(
          (x,y)->error("Don't know how to compute the newOp on ".expr2text(domtype(x))), 
                        Name="newOp");
\end{verbatim}
\normalsize and by assuming that the action on the Schur basis is
given by a function  \texttt{f}, we
declare
\scriptsize
\begin{verbatim}
 >>  operators::overloaded::declareSignature(
        S::newOp, [S::s, DOM_INT],
        S::s::moduleMorphism(f, dom::s));
\end{verbatim}
\normalsize and you can apply this operator on any basis due to
the overloading mechanism.

\section{Integration of others programs written in C++}

In July 2005, the American Institute of Mathematics organized in Palo
Alto, California, a workshop on the generalized Kostka polynomials
under the leadership of A. Schilling and M. Vazirani. During problem
sessions, MuPAD-Combinat was put to use for testing conjectures. As
usual in algebraic combinatorics, computations required the
combination of preexisting combinatorial and algebraic
functionnalities (as provided by MuPAD-Combinat) with new
combinatorial features (namely a highly technical bijection between
$k$-tuples of Young Tableaux and unrestricted rigged configurations).
Thanks to a dynamic module we could reuse a pre-existing robust C++
implementation of this bijection written by L. Deka; this allowed us
to start manipulating large examples quickly. In general, dynamic
modules permit us to reuse C++ code with two goals in mind: to avoid
reimplementing nontrivial and tested code, and to get quicker
computations than in pure MuPAD language. In section 4.1, we describe
the implementation of combinatorial objects in MuPAD, and in section
4.2 we present our technical choices for a seamless integration of
a combinatorial bijection implemented in C++.

\subsection{Implementation of combinatorial objects in MuPAD-Combinat}
We implement each combinatorial class as a domain in the MuPAD
category \texttt{Cat::CombinatorialClassWith2DBoxedRepresentation}.\\
This category provides, among other things, a pretty printing method
using ASCII characters. Let us illustrate this design in the case of
skew riggings implemented in the domain
\texttt{combinat::skewRiggings}. \scriptsize
\begin{verbatim}
  domain combinat::skewRiggings
      inherits Dom::BaseDomain;
      category Cat::CombinatorialClassWith2DBoxedRepresentation;
      axiom    Ax::canonicalRep;
    
      info_str := "Combinatorial class for rigged skew partitions";
  ...
  end_domain:
\end{verbatim}
\normalsize We implement the constructor 
\texttt{combinat::skewRiggings::new}  which builds an
object from the list of its operands:\scriptsize
\begin{verbatim}
  >> a:=combinat::skewRiggings([[], [[[0], [0, 1]], [[""], [0, 1]]]]);

                               +---+
                               | 0 | 0   1
                               +---+
                               |   | 0   1
                               +---+
\end{verbatim}
\normalsize The first element of the internal representation is the
type \scriptsize
\begin{verbatim}
  >> a[0];
  
                          combinat::skewRiggings

\end{verbatim}
\normalsize Ribbon rigged configurations are particular sequences of
skew riggings, implemented as plain lists of typed MuPAD objects. We
implement them in the following domain{\scriptsize
\begin{verbatim}
  domain combinat::riggedConfigurations::RcRibbonsTableaux
      inherits Dom::BaseDomain;
      category Cat::CombinatorialClassWith2DBoxedRepresentation, 
      // Elements of this domain are represented using a plain MuPAD data structure
               Cat::FacadeDomain(DOM_LIST);
    
      info_str := "Combinatorial class for rigged configurations";
  ...
  end_domain:
\end{verbatim}
}
Here is an example of the bijection applied on the set of
all $3$-ribbon tableaux of shape $(432)$ and evaluation $(111)$\scriptsize
\begin{verbatim}
  >> rc := map(combinat::ribbonsTableaux::list([4,3,2],[1,1,1],3), 
         combinat::riggedConfigurations::RcRibbonsTableaux::fromRibbonTableau);

  -- --                                   +---+     --
  |  |                                    |   |      |
  |  |  +---+          +---+---+          +---+---+  |
  |  |  | 1 | 0   1  , | 0 | 0 | 0   0  , |   |   |  |,
  |  -- +---+          +---+---+          +---+---+ --
  --

   --                                   +---+     --
   |                                    |   |      |
   |  +---+          +---+---+          +---+---+  |
   |  | 0 | 0   1  , | 0 | 0 | 0   0  , |   |   |  |,
   -- +---+          +---+---+          +---+---+ --

   --                               +---+     -- --
   |                 +---+          |   |      |  |
   |  +---+          | 0 | 0   0    +---+---+  |  |
   |  | 0 | 0   0  , +---+        , |   |   |  |  |
   |  +---+          |   | 0   1    +---+---+  |  |
   --                +---+                    -- --

  >> a:= rc[1];

             --                                   +---+     --
             |                                    |   |      |
             |  +---+          +---+---+          +---+---+  |
             |  | 1 | 0   1  , | 0 | 0 | 0   0  , |   |   |  |
             -- +---+          +---+---+          +---+---+ --
  >> op(a)[2][0];

                          combinat::skewRiggings

\end{verbatim}
\normalsize L. Deka and A. Schilling introduced in \cite{Schilling3} a
new kind of rigged configurations, namely the unrestricted ones. They
were implemented as an independent C++ program (file:  {\tt
  FromOneCrystalPath.cc}, headers:  {\tt
  FromOneCrystalPath.h}). On each rigged configurations
(ribbons one and unrestricted one) we can compute a statistic. The
interesting question is to find, in a special case, a bijection which
preserves the statistic between these two kinds of rigged
configurations. In order to manipulate these two objects in a single
program, we decided in collaboration with A.  Schilling and L.  Deka,
to also integrate this C++ program into MuPAD-Combinat. In order to
make the integration of this version of rigged configurations in a
transparent way for the user, we kept the same design we used for
ribbon rigged configurations.

\subsection{Implementation of combinatorial objects in a dynamic module using the MAPITL library}

We describe now the integration of the previous C++ program
using MAPITL library (MuPAD Application Programming Interface Template
Library (included in \texttt{MuPAD-Combinat}).  This library provides
\begin{itemize}
 \item Wrappers to use MuPAD lists as standard containers
 \item Easy C++ $\longleftrightarrow$ MuPAD conversions with one
   single overloaded template for each directions:
   \begin{itemize}
          \item C++ object $\longrightarrow$ MuPAD object: {\tt CtoM(c)}
          \item MuPAD object $\longrightarrow$ C++ object: {\tt MtoC(c)} 
          \end{itemize}
\end{itemize}
This includes conversions to/from containers and recursive containers.\\

The C++ program computes the rigged configurations corresponding to a
list of Young tableaux also called path which is denoted by the class
\texttt{path\_class} . We want to call from MuPAD \\ \texttt{void
  path\_class::build\_rigged\_for\_path} using the procedure\\
\texttt{combinat::RiggedConfigurations::newRiggedConfigurations}. \\The
building of the interaction is divided into three steps which are
realized in the file \\ \texttt{RiggedConfigurationsPaths.mcc}
\begin{itemize}
\item convert a list of Young Tableaux in MuPAD to a C++ object of the
  class  path\_class 
\item call the function  \texttt{void
    path\_class::build\_rigged\_for\_path}
\item convert the result into a MuPAD object of type 
  \texttt{DOM\_LIST} 
\end{itemize}
Original files are \texttt{FromOneCrystalPath.cc} containing
declarations of different classes. We now explain some parts of the
file \texttt{RiggedConfigurationsPaths.mcc} which is entirely given in
appendix. The beginning of the file is the inclusion of the header
file of the independent program and MAPITL \scriptsize
\begin{verbatim}
#include<vector>
#include "FromOneCrystalPath.h"
#include "MAPITL.h"
#include "MAPITL_tmpl.h"
\end{verbatim}
\normalsize
\subsubsection{Conversion MuPAD to C++}
\scriptsize
\begin{verbatim}
  MFUNC( newRiggedConfigurations, MCnop )
  {
    MFnargsCheck(3);
    MFargCheck(1, DOM_INT);
    MFargCheck(2, DOM_INT);
    MFargCheck(3, DOM_LIST);
    
    n = MtoC<int>(MFarg(1));
    int path_len  = MtoC<int>(MFarg(2));
    Cell input(MFarg(3));
    ...
    Cell *toto = input.toArray<Cell>(); 
    ...
}
\end{verbatim}
\normalsize  \texttt{MFargCheck}  is the type
checking of MuPAD arguments and  \texttt{MtoC} 
permits to convert MuPAD objects into C++ ones. After we initialize an
object of class  \texttt{path\_class}  with
values contained in  \texttt{toto}  \scriptsize
\begin{verbatim}
  path_class*  input_path = new path_class(path_len)
  for(long a=0; a < input.size(); a++)
    { 
      ...
    }
\end{verbatim}
\normalsize

\subsubsection{Using the original function in C++}
We can call now the original function in order to compute the
bijection
\scriptsize
\begin{verbatim}
  input_path->build_rigged_for_path()
\end{verbatim}
\normalsize
\subsubsection{Conversion C++ to MuPAD}
Next we create a C++ vector which contained operands of the
corresponding into MuPAD object 
\scriptsize
\begin{verbatim}
  for(i=0; i<n; i++)
   {
      vector<int> yyy;
      j = 0;
      while(input_path->rigged[i][0][j] != UNUSED)
      {
          yyy.push_back(input_path->rigged[i][0][j]);
          yyy.push_back(input_path->rigged[i][2][j]);
          yyy.push_back(input_path->rigged[i][1][j]);
          j++;
      }
      res[i] = yyy;
   }
\end{verbatim}
\normalsize
and we return the MuPAD list
\scriptsize
\begin{verbatim}
MFreturn(Cell(res));
\end{verbatim}
\normalsize
\subsection{Compilation of the dynamic module}
The next step is the compilation of the {\tt mcc} file giving us the
dynamic module \texttt{RiggedConfigurationsPaths.mdm}. We have to
load this module from MuPAD using the function \scriptsize
\begin{verbatim}
  combinat::RiggedConfigurationsPaths := 
  proc() 
    save RiggedConfigurationsPaths; 
  begin 
    if module::which("RiggedConfigurationsPaths") = FAIL then
      userinfo(1, "Dynamic module RiggedConfigurationsPaths not available"):
      RiggedConfigurationsPaths := FAIL;
    else 
      if traperror(module("RiggedConfigurationsPaths")) <> 0 then
         warning("Error loading the dynamic module RiggedConfigurationsPaths:");
         lasterror();
      end_if;
    end_if;
    RiggedConfigurationsPaths;
end_proc():
\end{verbatim}
\normalsize We create the domain \\
\texttt{combinat::riggedConfigurations::RcPathsEnergy }and procedure\\
\texttt{combinat::riggedConfigurations::RcPathsEnergy::fromOnePath }
computes the bijection by calling the function implemented in the dynamic module\\
\texttt{combinat::RiggedConfigurations::newRiggedConfigurations
}
\scriptsize
\begin{verbatim}
  >> a:=[combinat::tableaux([[3,3]]), combinat::tableaux([[2,2]]), combinat::tableaux([[1,1]])];

                   -- +---+---+  +---+---+  +---+---+ --
                   |  | 3 | 3 |, | 2 | 2 |, | 1 | 1 |  |
                   -- +---+---+  +---+---+  +---+---+ --
\end{verbatim}
\normalsize
We compute the bijection
\scriptsize
\begin{verbatim}
  >> rc := combinat::riggedConfigurations::RcPathsEnergy::fromOnePath(a);

                -- +---+---+                            --
                |  |   | 0 |     0    +---+---+          |
                |  +---+---+        , |   | 0 |     0    |
                |  |   | 0 |     0    +---+---+          |
                -- +---+---+                            --

  >> op(rc[1])[0]

                combinat::skewRiggings                                                 

\end{verbatim}
\normalsize The object computed using the C++ program is interfacing
in a transparent way for the user who can use other MuPAD
functionalities on the previous result.

\subsubsection{Acknowledgment}
The author wants to thank A. Schilling and L. Deka for their
collaboration on rigged configurations, M. Zabroki for its code on
creating operators, and F. Hivert and N. Thiery for their support in
the project MuPAD-Combinat.

\section*{Appendix: RiggedConfigurationsPaths.mcc}
\scriptsize
\begin{verbatim}
#include<vector>
#include "FromOneCrystalPath.h"
#include "MAPITL.h"
#include "MAPITL_tmpl.h"
using MAPITL::MtoC;
using MAPITL::CtoM;
using MAPITL::Container;
using MAPITL::SimpleCell;
using MAPITL::Cell;
using MAPITL::ObjectInCell;
using MAPITL::ArrayInCell;

using namespace std;

MFUNC( newRiggedConfigurations, MCnop )
{
  MFnargsCheck(3);
  MFargCheck(1, DOM_INT);
  MFargCheck(2, DOM_INT);
  MFargCheck(3, DOM_LIST);
  n = MtoC<int>(MFarg(1));
  int path_len  = MtoC<int>(MFarg(2));
  Cell input(MFarg(3));
  int i, j, k, tmp, path_index, tblu_index, col;
  
  tmp = UNUSED;
  path_index = 0;
  tblu_index = 0;
  i = 0; j = 0; k = 0; col = 0; l = 0;
  initialize_lambda();
  path_class*  input_path = new path_class(path_len);
  reset_tableau();
  Cell *toto = input.toArray<Cell>();  
  vector< vector<int> > res(n+1);
  
  for(long a=0; a < input.size(); a++)
  {  
      Cell* alpha=(toto[a]).toArray<Cell>();
      int s = (alpha[0]).size();
      tblu_class *my_tblu = new tblu_class((toto[a]).size() ,s);
      my_tblu->tblu_id = tblu_index;
      tblu_index += 1;
     
      for(long d=0; d<toto[a].size(); d++)   
      {
          int* beta=(alpha[d]).toArray<int>();
          for(long b=0; b < s ;b++)
          {
              my_tblu->tb[d][b] = beta[b];
              my_tblu->tab_lambda[beta[b]-1] =
                  my_tblu->tab_lambda[beta[b]-1] + 1; 
          }
      }
      input_path->path[path_index] = my_tblu;
      reset_tableau ();
      path_index += 1;
  }
  input_path->build_rigged_for_path();  
  for(i=0; i<n; i++)
  {
      vector<int> yyy;
      j = 0;
      while(input_path->rigged[i][0][j] != UNUSED)
      {
          yyy.push_back(input_path->rigged[i][0][j]);
          yyy.push_back(input_path->rigged[i][2][j]);
          yyy.push_back(input_path->rigged[i][1][j]);
          j++;
      }
      res[i] = yyy;
  }
  input_path->calculate_cocharge();
  vector<int> stat;
  stat.push_back(input_path->cocharge);
  res[n] = stat;
  
  MFreturn(Cell(res));
} MFEND
\end{verbatim}
\normalsize

\end{document}